\documentclass[12pt]{article}
\usepackage{amsmath}
\usepackage{amssymb}
\usepackage{amsthm}
\usepackage{graphicx}
\usepackage[usenames]{color}

\begin{document}
\vspace*{1.5cm}
\begin{center}
{\bf \Large A New Approach to Variational Inequalities\vspace{0.3cm}
of Parabolic Type}
\vspace{1cm}

Maria Gokieli, Nobuyuki 
Kenmochi and Marek Niezg\'odka \vspace{0.7cm}

  Interdisciplinary Centre for Mathematical
  and Computational Modelling,
  
   University of Warsaw, Pawi\'nskiego 5a,
    02-106 Warsaw, Poland 

\end{center}
\vspace{1cm}

\noindent
{\bf Abstract.} This paper is concerned with the
weak solvability of fully nonlinear parabolic
variational inequalities with time dependent
convex constraints. As possible approaches to such problems, there are for
instance the time-discretization method and the fixed point method of 
Schauder type with appropriate compactness theorems. In this paper, our attention is paid 
to the latter approach. However, there has not been prepared any appropriate 
compactness theorem up to date that enables us
the direct application of fixed point method to
variational inequalities of parabolic type. 
In order to establish it we have to start on the
set up of a new compactness theorem for a wide class of parabolic variational inequalities.
\vspace{1cm}

\noindent
{\bf 1. Introduction}
\vspace{0.5cm}

We consider a variational inequality of 
quasi-linear parabolic type:
  $$ \begin{array}{l}
  \displaystyle{\sum_{i=1}^2 \int_Q \frac{\partial u_i}
  {\partial t}  (u_i-\xi_i) dxdt 
  +\sum_{i=1}^2\sum_{k=1}^N 
  \int_Q a_i^{(k)}(x,t,u)\frac{\partial u_i} {\partial x_k} \frac{\partial (u_i-\xi_i)} {\partial x_k} dxdt }\\
   \displaystyle{~~~~~~~~~~~~~~~~~  \le \sum_{i=1}^2 \int_Q f_i(u_i-\xi_i)
  dxdt, }
    \end{array} \eqno{(1.1)} $$
  $$ \forall \xi:=[\xi_1.\xi_2]\in 
  L^2(0,T;H^1_0(\Omega)\times H^1_0(\Omega))
    ~{\rm with~} \xi(t) \in K(t) ~{\rm a.e.~}t
         \in [0,T], 
   \eqno{(1.2)}$$
  $$ u(x,0) =u_0(x)~~
  {\rm in~}\Omega,~~u=0~{\rm on~}\Sigma,
     \eqno{(1.3)}$$
where $\Omega$ is a bounded domain in 
${\bf R}^N$, $Q:=\Omega\times (0,T),
~0<T<\infty$, $\Gamma:=\partial \Omega$, 
$\Sigma :=\Gamma \times (0,T)$, $u:=[u_1,u_2]$ and the diffusion coefficients 
$a_i^{(k)}(x,t,u)$ are strictly positive, bounded and continuous in $(x,t,u)\in \overline Q\times {\bf R}^2$ as well as constraint set
$K(t)$ is convex and closed in 
$H^1_0(\Omega)\times H^1_0(\Omega)$
satisfying some smoothness assumption in 
$ t \in [0,T]$.
Functions $f:=[f_1,f_2]$ and $u_0$ are
 prescribed 
in $L^2(Q)\times L^2(Q)$ and $K(0)$,
 respectively, as the data.
Our claim is to construct a solution $u$ of (1.1)-(1.3) in a weak sense such that
 $$ u \in C([0,T];L^2(\Omega)\times L^2(\Omega))\cap L^2(0,T;H^1_0(\Omega)\times H^1_0(\Omega)),~~u(t) \in K(t), ~\forall
   t \in [0,T].$$

In the case without constraint, namely $K(t)=H^1_0(\Omega)\times H^1_0(\Omega)$, our
problem is the usual initial-boundary value problem for parabolic system of quasi-linear
PDEs :
  $$ \frac{\partial u_i}{\partial t}
   -\sum_{k=1}^N \frac{\partial }{\partial x_k}\left(a_i^{(k)}(x,t,u)
    \frac{\partial u_i}{\partial x_k}\right) =f_i(x,t) ~~{\rm in~}Q,~i=1,2. $$ 
For the solvability a huge number of results have been established 
(cf. [14, 24]), for instance, the Leray-Schauder principle together with some 
compactness theorems, such as [3, 26, 31].

In connection with quasi-linear variational
inequalities, the concept of nonlinear monotone
mappings was generalized to several classes of
nonlinear mappings of monotone type, for
instance, semimonotone
[25], pseudomonotone 
[5, 11, 17], and furthermore $L-$pseudomonotone
mappings [6, 7]. Especially the last one was
introduced as a class of nonlinear perturbations for linear maximal monotone 
mappings $L$, which
is available for parabolic variational inequalities; in a typical application 
of this theory, $L$ is the time-derivative
$\frac d{dt}$.   
However, it seems still difficult
to treat directly our model problem (1.1)-(1.3)
in these frameworks of nonlinear mappings of monotone type. \vspace{0.5cm} 

Our model problem (1.1)-(1.3) is formally written in the space 
$L^2(0,T;X^*)$ with $X:=L^2(0,T;H^1_0(\Omega)\times H^1_0(\Omega))$ ($X^*=L^2(0,T;H^{-1}(\Omega)\times H^{-1}(\Omega))$) as
    $$f \in Lu+A(u,u),~~u(0)=u_0, $$
by taking as $L$ the mapping $L:=\frac d{dt}+\partial I_{K(t)}(\cdot)
: D(L) \subset X \to X^*$ 
and as $A$ the mapping $A(v,u): D(A)=X \to X^*$
given by
 $$ \langle A(v,u),[\xi_1,\xi_2]
 \rangle _{X^*,X}= 
  \sum_{i=1}^2\sum_{k=1}^N \int_Q a_i^{(k)}(x,t,v)\frac {\partial u_i}{\partial x_k}
\frac {\partial \xi_i}{\partial x_k} dxdt, $$
   $${\rm for}~ u:=[u_1,u_2],~v=[v_1,v_2],
   ~\xi=[\xi_1,\xi_2] \in X, $$
where $\langle \cdot, \cdot\rangle_{X^*,X}$ stands for the duality between $X^*$ and $X$.
We see that $L$ is maximal monotone from $D(L) \subset X$ into
$X^*$, but $L$ is nonlinear in general.
Since 1970, it remains to set up an abstract approach
to such a quasi-linear parabolic variational inequality as our model problem. 
In this paper we establish a new approach 
to parabolic variational inequalities
with time-dependent constraints $\{K(t)\}$, based on a new compactness theorem 
(see Theorem 2.1) derived from the total variation estimates for solutions
of parabolic variational inequalities; this idea was found in a recent work [16] of the authors. 

There is a different approach to nonlinear variational inequalities of 
parabolic type with time-independent convex constraint in [1] where 
the time-discretization method was employed and a compactness theorem
 was established to ensure the strong convergence of 
time-discretized approximation schemes in time. This idea seems available to 
the case of time-dependent convex constraints.

In this paper the following notations are used. For a general (real) Banach space $X$ we denote
by $X^*$ the dual space, by $|\cdot|_X$ and
$|\cdot|_{X^*}$ the norms in $X$ and $X^*$,
respectively, and by $\langle \cdot, \cdot \rangle_{X^*,X}$ the duality between $X^*$ and $X$. Especially, when $X$ is a Hilbert space,
we denote by $(\cdot,\cdot)_X$ the inner product
in $X$.

Let $\varphi(\cdot)$ be a proper, lower 
semi-continuous (l.s.c.) and convex function on
a Banach space $X$. Then the subdifferential
$\partial_*\varphi$ of $\varphi$ is defined by
a multivalued mapping from $X$ into $X^*$ as
follows: $z^* \in \partial_*\varphi(z)$ if and
only if $z \in D(\varphi):=\{z \in X~|~\varphi
(z)<\infty\}$, $z^* \in X^*$ and
 $$ \langle z^*, w-z\rangle_{X^*,X} \le 
     \varphi(w) -\varphi(z),~~\forall w\in X.$$
The set $D(\varphi)$ is called the effective
domain of $\varphi$. The domain of
$\partial_*\varphi$ is the set 
$D(\partial_*\varphi):=\{z \in D(\varphi)~|~
\partial_*\varphi(z)\ne \emptyset\}$ and the
range of $\partial_*\varphi$ is the set
$R(\partial_*\varphi):=
\cup_{z \in D(\varphi)}
\partial_*\varphi(z)$.

Let $K$ be a non-empty closed and convex subset
 of a Banach space $X$. Then the function $I_K(\cdot)$ given by
 $$  I_K(z):=\left \{
      \begin{array}{l}
       0,~~~~{\rm if~}z \in K,\\[0.3cm]
        \infty,~~~{\rm if~}z \in X-K,
      \end{array}\right. $$
is proper, l.s.c. and convex on $X$ and is
called the indicator function of $K$ on $X$,
and the subdifferential $\partial_*I_K$ is
defined as a multivalued mapping from $X$ into
$X^*$. Clearly $D(I_K)=D(\partial_*I_K)=K$.

For the real function $r \to |r|^{p-1}$ on
${\bf R}_+:=\{r \in {\bf R}~|~r \geq 0\}$ for a
fixed number $p > 1$ and
a Banach space $X$, we
consider the mapping $F: X \to X^*$ which
assigns to each $z \in X$ the set
 $$Fz:=\{z^*\in X^*~|~\langle z^*,z\rangle
  _{X^*,X}=|z|^p_X, |z^*|_{X^*}=|z|_X^{p-1}\}.$$
This mapping $F$ is called the duality mapping
associated with the gauge function 
$r \to |r|^{p-1}$. It is well known that $F$ is the subdifferential of the 
non-negative, continuous and convex
function $\psi(z):= \frac 1p |z|^p_X$ on $X$,
namely, $Fz=\partial_*\psi(z)$ for all 
$z \in D(F)=X$. In particular, if $X$ is reflexive
and $X^*$ is strictly convex, then $F$ is singlevalued and continuous 
from $X$ into $X^*_w$ ($X^*_w$ stands for the space $X^*$ 
with the weak topology); this continuity is called the ^^ ^^ demicontinuity" of $F$. 
Let $A$ be a (multivalued) mapping from a Banach
space $X$ into $X^*$; the graph $G(A)$ of
$A$ is the set
 $$ G(A):=\{[z,z^*] \in X\times X^*~|~z^*\in Az\}.$$
Then $A$ is called monotone from $D(A)\subset X$ into $X^*$, 
 if
 $$ \langle z^*_1-z^*_2, z_1-z_2\rangle \geq 0,
   ~~\forall [z_i,z^*_i]\in G(A),~i=1,2;$$
in particular, $A$ is called strictly monotone,
if the strict positiveness holds whenever 
$z_1 \ne z_2$ in the above inequality.
Moreover, $A$ is maximal monotone, if $A$ is
monotone from $D(A)\subset X$ into $X^*$ and 
$G(A)$ has no proper monotone extension in
$X\times X^*$. When $X$ is reflexive, it is well known that $A$ is maximal
monotone if and
only if the range of $A+F$ is the whole of 
$X^*$, where $F$ is the duality mapping from
$X$ into $X^*$. We refer to [2, 4, 5, 9, 10, 20, 26, 27] for 
fundamental properties of subdifferentials and monotone mappings.
\vspace{1cm}

\noindent
{\bf 2. A compactness theorem}\vspace{0.5cm}

In this section, let $V$ be a (real) reflexive
Banach space which is dense and compactly
embedded in a Hilbert space $H$. Identifying 
$H$ with its dual space, we have
$ V\subset H\subset V^*$ with compact embeddings.
Let $W$ be another reflexive and separable
 Banach space which
is dense and continuously embedded in $V$; 
since $V^* \subset W^*$, it holds that
   $$ V\subset H\subset W^*~~{\rm with~compact
 ~embeddings}. $$

In order to avoid
some irrelevant arguments, suppose that 
$V,~V^*,W$ and $W^*$ are strictly convex.
We denote 
by $C_W$ an embedding constant from $W$ into 
$V$ and $H$, namely 
 $$|z|_V \le C_W|z|_W,~~|z|_H \le C_W|z|_W, 
    ~~\forall z \in W.$$

For any function $w:[0,T]\to W^*$,
the total variation of $w$ is denoted by 
${\rm Var}_{W^*}(w)$,
which is defined by
 $$ {\rm Var}_{W^*}(w):=\sup_{
    \begin{array}{l}
    \eta \in C^1_0(0,T;W),\\
   |\eta|_{L^\infty(0,T;W)}\le 1
     \end{array}} \int_0^T
  \langle w,\eta'\rangle_{W^*,W}dt. $$
We refer to [9] or [13] for the fundamental properties of total
variation functions.\vspace{0.5cm}

We fix numbers $p$ with $1<p<\infty$, 
$p':=\frac p{p-1}$, and 
$T$ with $0<T<\infty$.
Given $\kappa>0,~M_0 >0$ and $u_0\in H$, consider the
set $Z(\kappa, M_0,u_0)$ in $L^p(0,T;V)\cap 
L^\infty(0,T;H)$ given by:
$$
  Z(\kappa, M_0,u_0):=\left\{u~\left|
   \begin{array}{l}
     |u|_{L^p(0,T;V)}\le M_0,~
 |u|_{L^\infty(0,T;H)}\le M_0,\\[0.2cm]
  \exists f \in L^{p'}(0,T;V^*)~{\rm such~that~}
  \\[0.2cm]
 \displaystyle{~~~\int_0^T \langle f, 
   u\rangle dt
  \le M_0,~|f|_{L^1(0,T;W^*)} \le M_0,}\\[0.3cm]
  \displaystyle{~~~ \int_0^T \langle \eta'-f, u-\eta 
 \rangle dt\le 
   \frac 12 |u_0-\eta(0)|^2_H,~\forall \eta 
   \in L^p(0,T;V)}\\[0.2cm]
   ~~~~~~
   {\rm with~}\eta' \in L^{p'}(0,T;V^*),~
  \eta(t) \in \kappa B_W(0),~\forall t \in[0,T]
  \end{array} \right. \right\},
 $$
where $B_W(0)$ is the closed unit ball in $W$
with center at the origin and $\langle \cdot, \cdot \rangle=
\langle \cdot, \cdot \rangle_{V^*,V}$. \vspace{0.5cm}

The variational inequality
 $$\int_0^T \langle \eta'-f, u-\eta 
 \rangle dt\le 
   \frac 12 |u_0-\eta(0)|^2_H, \eqno{(2.1)} $$
 $$\forall \eta \in L^p(0,T;V),~\eta' \in 
 L^{p'}(0,T;V^*),~\eta(t) \in \kappa B_W(0),
 ~\forall t \in[0,T], \eqno{(2.2)}$$
in the definition of $Z(\kappa,M_0,u_0)$,
is derived from the time-derivative with the
convex constraint $\kappa B_W(0)$. In fact, for
any $u,~\eta \in L^p(0,T;V)$ with $u',~\eta' 
\in L^{p'}(0,T;V^*)$ with $u(0)=u_0$ 
we have by integration by parts
 $$ \int_0^T\langle \eta' -u', u-\eta \rangle dt
 =\frac 12|u(0)-\eta(0)|^2_H
   -\frac 12|u(T)-\eta(T)|^2_H
   \le \frac 12|u_0-\eta(0)|^2_H.$$
Therefore, if $f=u'$ with $u(0)=u_0$, then (2.1) holds. 
Given $u \in L^p(0,T;V)\cap L^\infty(0,T;H)$ and
$u_0 \in H$, the set of all $f$ satisfying (2.1)-(2.2) includes 
$u'$, provided $u'$ exists in $L^{p'}(0,T;V^*)$ and
$u(0)=u_0$. However, in general, it is an
extremely
large set; note that in the definition of
$Z(\kappa,M_0,u_0)$, any differentiability of
$u$ in time is not required.\vspace{0.5cm}

Our main result is stated as follows. \vspace{0.5cm}

\noindent
{\bf Theorem 2.1.} 
{\it Let $\kappa>0,~M_0 >0$ be any numbers and $u_0$ be any element of $H$. Then the set
$Z(\kappa,M_0,u_0)$ is relatively compact in
 $L^p(0,T;H)$. Moreover, the convex closure of $Z(\kappa,M_0,u_0)$,
denoted by ${\rm conv}[Z(\kappa, M_0,u_0)]$, in $L^p(0,T;V)$ is compact
in $L^p(0,T;H)$.} \vspace{0.5cm}

We begin with the
following lemmas that are crucial for the proof 
of Theorem 2.1.\vspace{0.5cm}

\noindent
{\bf Lemma 2.1.} {\it Let $\kappa>0,~M_0 >0$ be
any numbers and $u_0$ be any element of $H$.
Then there exists a positive constant 
$C^*:=C^*(\kappa,M_0, |u_0|_H)$, 
depending only on 
$\kappa, ~M_0$ and $|u_0|_H$, such that
  $$ {\rm Var}_{W^*}(u) \le C^*, \eqno{(2.3)}$$
for all $u \in Z(\kappa, M_0, u_0)$. }
\hspace*{-0.3cm}

\noindent
{\bf Proof.} Let $u$ be any element in 
$Z(\kappa, M_0, u_0)$,
and take a function $f \in L^{p'}(0,T;V^*)$
satisfying all the required properties in the
definition of $Z(\kappa, M_0,u_0)$.
Now let $\eta$ be any function in $C^1_0(0,T;W)$
with $|\eta|_{L^\infty(0,T;W)}>0$.
Since $\tilde \eta :=\pm \frac {\kappa \eta}
{|\eta|_{L^\infty(0,T;W)}}$ is a possible test
function for (2.1)-(2.2), we have
 $$ \int_0^T \langle \tilde \eta'-f, u-
 \tilde \eta \rangle dt
   \le \frac 12 |u_0|^2_H,$$
which shows
 \begin{eqnarray*}
 \int_0^T\langle \tilde \eta', u \rangle dt
&\le& -\int_0^T\langle f,\tilde\eta \rangle dt
 +\int_0^T\langle f, u \rangle dt
  +\frac 12 |u_0|^2_H,\\
&\le& |f|_{L^1(0,T;W^*)}
  |\tilde \eta|_{L^\infty(0,T;W)}+M_0
  +\frac 12 |u_0|^2_H.
  \end{eqnarray*}
Hence,
 $$ \left |\int_0^T\langle u, \eta' \rangle dt
    \right|
    \le \left\{M_0+
   \frac {M_0}\kappa+\frac {|u_0|^2_H}
   {2\kappa}\right\}
    |\eta|_{L^\infty(0,T;W)},~~
  \forall \eta \in C^1(0,T;W),$$
so that (2.3) holds with 
$C^*:= M_0
  +\frac 1{\kappa}M_0
   +\frac 1{2\kappa} |u_0|^2_H$.   \hfill $\Box$
\vspace{0.5cm}

\noindent
{\bf Lemma 2.2.} {\it Let $M_1$ be any positive number 
and let $\{u_n\}$ be any sequence of functions from $[0,T]$ into $W^*$
such that $u_n \in L^p(0,T;V)\cap L^\infty(0,T;H)$
 $$ |u_n|_{L^p(0,T;V)} \le M_1,~~|u_n|_{L^\infty(0,T;H)}\le M_1,
      ~~{\rm Var}_{W^*}(u_n) \le M_1,~~n=1,2,\cdots. \eqno{(2.4)}$$ 
Then there are a subsequence $\{u_{n_k}\}$ of
$\{u_n\}$ and a function $u \in L^p(0,T;V)\cap L^\infty(0,T;H)$
such that $u_{n_k}(t) \to u(t)$ weakly in $H$
for every $t \in [0,T]$ as $k \to \infty$.
Hence $u_{n_k}(t) \to u(t)$ in $W^*$
for every $t \in [0,T]$ 
and $u_{n_k} \to u$ in $L^q(0,T;W^*)$ for
 every $q \in [1,\infty)$ as $k \to \infty$.}\vspace{0.1cm}

\noindent
{\bf Proof.} 
Since $W$ is
separable, there is a countable dense subset 
$W_0$ of $W$. 
Now, we consider a 
sequence of real valued functions $A_n(t,\xi):=(u_n(t), \xi)_H~
(=\langle u_n(t),
\xi \rangle_{W^*,W})$ on $[0,T]$ for each 
$\xi \in W_0$.
Then, by (2.4) the total
variation of $A_n(t,\xi)$ is bounded by 
$M_1|\xi|_W$. 
Hence from the Helly
selection theorem (cf. [13 ; Section 5.2.3]) 
it follows that there is a 
subsequence $\{n_k\}$, depending on 
$\xi \in W_0$, such that $A_{n_k}(t,\xi)$
converges to a function $A_0(t,\xi)$ pointwise on $[0,T]$ and its total 
variation is not larger than $M_1|\xi|_W$. 

Since $W_0$ is 
countable in $W$, by using extensively the above
Helly selection theorem
we can extract a subsequence, denoted by
the same notation as $\{n_k\}$ again, and
a function $A_0(t,\xi)$ on $[0,T]\times W_0$
such that
 $$ A_{n_k}(t,\xi) \to A_0(t, \xi)~{\rm as~}
 k\to \infty, ~~
   \forall t \in [0,T],~\forall \xi \in W_0.
 \eqno{(2.5)}$$ 
Furthermore, by density, this convergence (2.5)
can be extended to all $\xi \in W$. Also,
the functional $A_{n_k}(t,\xi)$ is linear in 
$\xi$ and uniformly bounded, i.e.
 $$ |A_{n_k}(t,\xi)| \le M_1|\xi|_H 
   \le M_1C_W|\xi|_W,~~\forall t \in [0,T],~
   \forall \xi \in W.$$
This implies that $A_0(t,\xi)$ is linear and bounded in $\xi \in W$ and 
$|A_0(t,\xi)| \le M_1|\xi|_H$ for all $\xi\in W$
and $t \in [0,T]$. As a consequence, by
Riesz representation theorem, there is a function 
$u: [0,T]\to H$ with $| u(t)|_H
\le M_1$ for all $t \in [0,T]$ such that
 $$ A_0(t,\xi)=(u(t), \xi)_H,~~\forall \xi\in
   H,~\forall t \in [0,T].$$
Now it is clear by (2.5) 
that $u_{n_k}(t) \to u(t)$ weakly in $H$ for
$t \in [0,T]$ as $k \to \infty$. Finally, by the compactness of the
injection from $H$ into $W^*$, we see that $u_{n_k}(t) \to u(t)$ (strongly) in $W^*$ for
any $t \in [0,T]$. Hence $u_{n_k} \to u$ in $L^q(0,T;W^*)$ for all
$q \in [1,\infty)$ as $k \to \infty$.
\hfill $\Box$
\vspace{0.5cm}

\noindent
{\bf Proof of Theorem 2.1.} We first note from Lemma 2.1 that
  $$ Z(\kappa,M_0,u_0) \subset {\cal X}:=
     \{u~|~|u|_{L^p(0,T;V)}\le M_0,~|u|_{L^\infty(0,T;H)}\le M_0,~
     {\rm Var}_{W^*}(u) \le C^* \},$$
where $C^*$ is the same constant as in Lemma 2.1. Note that ${\cal X}$ is 
closed and convex in $L^p(0,T;V)$.
Therefore, in order to obtain Theorem 2.1 it is
enough to prove the compactness of ${\cal X}$ in $L^p(0,T;H)$.

Let $\{u_n\}$ be any
sequence in the set ${\cal X}$.
Then, by Lemma 2.2, there is a subsequence
$\{u_{n_k}\}$ and a function $u \in L^\infty
(0,T;H)$ such that $u_{n_k}(t) \to u(t)$ weakly
in $H$ for every $t\in [0,T]$ as $k \to \infty$.
By the injection compactness from $H$ into
$W^*$ we have that
 $$u_{n_k} \to u~~{\rm in~}L^p(0,T;W^*)~{\rm
   as~}k\to \infty. \eqno{(2.6)}$$
and that $|u|_{L^p(0,T;V)} \le M_0$ by
$|u_{n_k}|_{L^p(0,T;V)} \le M_0$.

Here we recall the Aubin lemma [3] (or [26; Lemma 5.1]):
for each $\delta>0$ there is a positive constant
$C_\delta$ such that
 $$ |z|^p_H \le \delta |z|^p_V 
  + C_\delta |z|^p_{W^*},
    ~~\forall z \in V.$$
By making use of this inequality for 
$z=u_{n_k}(t)-u(t)$ and integrating it in time, we get
 $$ \int_0^T |u_{n_k}(t)-u(t)|^p_H dt  \le 
   \delta (2M_0)^p
 + C_\delta \int_0^T|u_{n_k}(t)-u(t)|^p_{W^*}
   dt.$$
On account of (2.6), letting 
$k \to \infty$ gives that
$$\limsup_{k\to \infty}\int_0^T |u_{n_k}-u|^p_H dt
    \le \delta (2M_0)^p.$$
Since $\delta>0$ is arbitrary, we conclude that
$u_{n_k} \to u$ in $L^p(0,T;H)$. \hfill $\Box$
\vspace{0.5cm}

\noindent
{\bf Remark 2.1.} If $f \in L^{p'}(0,T;V^*)$
and (2.1) holds for all $\eta \in L^p(0,T;V)$ with 
$\eta' \in L^{p'}(0,T; V^*)$, then 
 $u'=f \in L^{p'}(0,T;V^*)$ and $u(0)=u_0$. Therefore,
by Theorem 2.1 the set
$$ \{u~|~
   |u|_{L^p(0,T;V)} \le M_0,~|u|_{L^\infty(0,T;H)}\le M_0,~
   |u'|_{L^1(0,T;W^*)} \le M_0 \}
  $$
is relatively compact in $L^p(0 ,T;H)$ for each
finite positive constant $M_0$. Our theorem includes
a typical case of Aubin compactness theorem 
[3]. \vspace{0.5cm}

\noindent
{\bf Remark 2.2.} A compactness result
of the Aubin type was extended in various directions, for instance
[12] and [18], and further to a quite general set up [31].
 \vspace{1cm}

\noindent
{\bf 3. Time-derivative under convex constraints} \vspace{0.5cm}

Let $H$ be a Hilbert space and
$V$ be a strictly convex reflexive Banach space
such that $V$ is dense in $H$ and the injection
from $V$ into $H$ is continuous. We identify $H$ with its dual space:
 $$ V\subset H \subset V^*~~{\rm with~
 continuous~ embeddings}.$$
For simplicity, we assume that $V^*$ is strictly convex. 
Therefore
the duality mapping $F$ from $V$ into $V^*$, associated with gauge function 
$r \to  |r|^{p-1}$, is singlevalued and
demicontinuous from $V$ into $V^*$, where $p$
is a fixed number with $1 < p <\infty$.

For the sake of simplicity for notation, we
write $\langle\cdot,\cdot\rangle$ for
$\langle\cdot,\cdot\rangle_{V^*,V}$ again.
\vspace{0.5cm}

Let $\{K(t)\}_{t \in [0,T]}$ be a family of non-empty, closed and
convex sets in $V$ such that there are functions 
$\alpha \in W^{1,2}(0,T)$ and
$\beta \in W^{1,1}(0,T)$ satisfying the following property:
for any $s, t \in [0,T]$ and any
$z \in K(s)$ there is $\tilde z \in K(t)$ such that
 $$ |\tilde z -z|_H \le |\alpha(t)-\alpha(s)|
 (1+|z|^{\frac p2}_V),~~
 |\tilde z|^p_V - |z|^p_V \le |\beta(t)-\beta(s)|
 (1+|z|^p_V).\eqno{(3.1)}$$

We denote by $\Phi(\alpha,\beta)$ the set of all such families $\{K(t)\}$, and put
 $$\Phi_S:=\bigcup_{\alpha \in W^{1,2}(0,T),~
 \beta\in W^{1,1}(0,T)} \Phi(\alpha,\beta),
  $$
which is called the strong class of 
time-dependent convex sets. \vspace{0.5cm}

Given $\{K(t)\} \in \Phi_S$, we consider the following time-dependent 
convex function on $H$:
  $$ \varphi_K^t(z):= \left \{
        \begin{array}{l}
 \displaystyle {\frac 1p |z|^p_V+I_{K(t)}(z), ~~~{\rm if~}z \in K(t),}\\[0.2cm]
 \displaystyle {\infty,~~~~~~~~~~~~~~~~~~~~~{\rm otherwise},}
         \end{array} \right. \eqno{(3.2)}$$
where $I_{K(t)}(\cdot)$ is the indicator function of
$K(t)$ on $H$. For each $t \in [0,T]$, 
$\varphi^t_K(\cdot)$
is proper, l.s.c. and strictly convex on $H$ and on $V$. 
By the general theory on nonlinear
evolution equations generated by time-dependent 
subdifferentials, condition (3.1) is 
sufficient in order
that for any $u_0 \in \overline{K(0)}$ (the closure of $K(0)$ in $H$) and $f \in L^2(0,T;H)$ the Cauchy problem
 $$ u'(t)+\partial \varphi^t_K(u(t))
 \ni f(t),~~u(0)=u_0,~{\rm in~}H,$$
admits a unique solution $u$ such that
$u \in C([0,T];H)\cap L^p(0,T;V)$ with $u(0)=u_0$, 
$t^{\frac 12}u'\in L^2(0,T;H)$ and 
$t \to t \varphi^t_K(u(t))$ is bounded on 
$(0,T]$ and absolutely continuous on any
 compact interval in $(0,T]$,
where $\partial \varphi^t_K$ denotes the subdifferential of 
$\varphi^t_K$ in $H$. In particular, 
if $u_0 \in K(0)$,
then $u' \in L^2(0,T;H)$ and $t \to \varphi^t_K(u(t))$ is
absolutely continuous on $[0,T]$.  \vspace{0.5cm}

Next, taking constraints of obstacle type into account, we introduce a weak 
class of time-dependent convex sets. In the sequel,
let $1<p<\infty$ and $\frac 1p+\frac 1{p'}=1$. 
\vspace{0.5cm}

\noindent
{\bf Definition 3.1.} Let $c_0$
be a fixed constant and $\sigma_0$ be a fixed 
function in $C([0,T];V)$ with $\sigma'_0 \in L^{p'}(0,T;V^*)$. Associated
with these $c_0$ and $\sigma_0$,  for each
small positive number $\varepsilon$
a mapping ${\cal F}_\varepsilon: [0,T]\times
V \to V$ is defined by
 $$ {\cal F}_\varepsilon(t)z 
=(1+ \varepsilon c_0)z +\varepsilon \sigma_0(t),
  ~\forall t \in [0,T],~\forall z \in V.
  \eqno{(3.3)}$$ 
Then, with ${\cal F}_\varepsilon$
the weak class $\Phi_W:=\Phi_W(c_0,\sigma_0)$ of time-dependent
convex sets is defined by:
$\{K(t)\}\in \Phi_W$ if and only if 
\begin{description}
\item {(a)} $K(t)$ is a closed and convex set in $V$ for all 
$t \in[0,T]$,
\item {(b)} there exists a
sequence $\{\{K_n(t)\}\}_{n \in {\bf N}} \subset \Phi_S $ such that for any 
$\varepsilon \in (0,\varepsilon_0]~(0<\varepsilon_0 <1) $ 
there is a positive integer 
$N_\varepsilon$ satisfying
 $$ {\cal F}_\varepsilon(t) (K_n(t)) \subset K(t),~~
{\cal F}_\varepsilon(t) (K(t)) \subset K_n(t),~~ 
      \forall t \in [0,T],~~\forall n\geq 
   N_\varepsilon.$$
\end{description}
In this case, it is said that $\{K_n(t)\}$
converges to $\{K(t)\}$ as $n \to\infty$,
which is denoted by 
 $$ K_n(t)\Longrightarrow K(t)~~{\rm on~}[0,T]~({\rm as~}n \to \infty).
   $$ 
\vspace{0.1cm}

We give three typical examples of $\{K(t)\}$ in the weak class $\Phi_W$.
\vspace{0.3cm}

\noindent
{\bf Example 3.1.} Let $\Omega$ be a bounded smooth domain in 
${\bf R}^N,~1 \le N <\infty$ and $Q:=\Omega \times (0,T)$. Let
$H:=L^2(\Omega),~V:=W^{1,p}(\Omega),~2 \le p<\infty$. Moreover, let 
$\rho:=\rho(x,t) \in C(\overline Q)$ and choose a sequence $\{\rho_n\}$ in 
$C^2(\overline Q)$ such 
that $\rho_n \to \rho$ in  $C(\overline Q)$. Now, constraint sets
$K(t)$ and $K_n(t)$ are defined by
 $$ K(t):=\{z \in V~|~z(x) \geq \rho(x,t)~{\rm for~a.e.~}x \in \Omega\},~
   \forall t \in [0,T],$$
and 
$$ K_n(t):=\{z \in V~|~z(x) \geq \rho_n(x,t)~{\rm for~a.e.~}x \in \Omega\},~
   \forall t \in [0,T].$$
Given $\varepsilon >0$, take a positive integer
$N_\varepsilon$ so that
    $$ |\rho_n-\rho |\le \varepsilon 
     ~~{\rm on~}\overline Q,~\forall n \geq 
     N_\varepsilon. $$
In this case, with the choice of $c_0=0$ and $\sigma\equiv 1$ the mapping 
${\cal F}_\varepsilon(t)$ is of the form
$z \to z+\varepsilon$, which maps $V$ into itself. Then we have:

(i) $\{K_n(t)\} \in \Phi_S$. Indeed,
for $z \in K_n(s)$, the function $\tilde z(x):=
z(x)-\rho_n(x,s)+\rho_n(x,t)$ belongs to 
$K_n(t)$ and (3.1) holds with functions 
$\alpha(t)=\beta(t)=c_nt$ for a constant $c_n>0$,
depending only on $|\rho_n|_{C^2(\overline Q)}$.
Thus $\{K_n(t)\} \in \Phi_S$.

(ii) $\{K(t)\} \in \Phi_W$. 
In fact, for any $z_n \geq \rho_n(\cdot,t)$ a.e. on
$\Omega$, we have
 $$ {\cal F}_\varepsilon z_n =z_n+\varepsilon
   \geq \rho_n(\cdot,t)+\varepsilon \geq \rho(\cdot,t)
   ~{\rm a.e.~on~} \Omega,$$
which implies ${\cal F}_\varepsilon (K_n(t))\subset K(t)$. Similarly, 
${\cal F}_\varepsilon (K(t)) \subset K_n(t)$. Hence 
$K_n(t) \Longrightarrow K(t)$ on $[0,T]$, and thus
$\{K(t)\} \in \Phi_W$. 
\vspace{0.5cm} 

\noindent
{\bf Example 3.2.} Let $\rho$ and $\rho_n$ be
the same as in Example 3.1, and consider
constraint sets
 $$ K(t):=\{z \in V~|~z(x) \le \rho(x,t)~{\rm for~a.e.~}x \in \Omega\},~
   \forall t \in [0,T]$$
and 
$$ K_n(t):=\{z \in V~|~z(x) \le \rho_n(x,t)~{\rm for~a.e.~}x \in \Omega\},~
   \forall t \in [0,T].$$
Then, just as in Example 3.1, we have
$\{K_n(t)\} \in \Phi_S$ and 
$K_n(t) \Longrightarrow K(t)$ on $[0,T]$
by using the mapping 
${\cal F}_\varepsilon(t)z=z-\varepsilon$, 
so that $\{K(t)\} \in \Phi_W$.
\vspace{0.5cm}

\noindent
{\bf Example 3.3.} Let $\Omega$ and $Q$ be the same as in Example 3.1 and
consider the following vectorial case in connection with our model problem
(1.1)-(1.3):
 $$H:=L^2(\Omega)\times L^2(\Omega),~
    V:=H^1_0(\Omega)\times H^1_0(\Omega);$$
hence $V^*=H^{-1}(\Omega)\times H^{-1}(\Omega)$.

Let $\psi=\psi(x,t)$ be an obstacle function prescribed in $C(\overline Q)$ 
so that $\psi \geq c_\psi$ on $\overline Q$ for a positive constant 
$c_\psi$, and define $K(t)$ by
 $$ K(t):=\{[\xi_1,\xi_2] \in V~|~
 |\xi_1|+|\xi_2| \le \psi(\cdot,t)~{\rm a.e.~in~}
 \Omega\},~~\forall t\in [0,T].$$
Next, choosing a sequence $\{\psi_n\}$ in 
$C^1(\overline Q)$ such that
 $$  \psi_n \geq c_\psi ~{\rm on~}\overline Q,
 ~\psi_n \to \psi ~{\rm in~}C(\overline Q)~~
 {\rm as~}
 n \to \infty$$
we define
 $$ K_n(t):=\{[\xi_1,\xi_2] \in V~|~
 |\xi_1|+|\xi_2| \le \psi_n(\cdot,t)~{\rm a.e.~in~}
 \Omega\},~~\forall t\in [0,T].$$
and the mapping 
 $${\cal F}_\varepsilon(t)[\xi_1,\xi_2]:=
  (1-\varepsilon)[\xi_1,\xi_2],~\forall 
  ~{\rm small~} \varepsilon >0;$$
note that this mapping is obtained by
the choice of $c_0=-1$ and $\sigma_0\equiv 0$.
Then we have:

(i) $\{K_n(t)\} \in \Phi_S$. In fact, given 
$z=[\xi_1,\xi_2] \in K_n(s)$, we take $\tilde z=[\tilde\xi_1,\tilde \xi_2]
:=(1-\frac 1{c_\psi}|\psi_n(s)-\psi_n(t)|_{C(\overline \Omega)})z$. 
In this case, if $|s-t|$ is so small that 
$\frac 1{c_\psi}|\psi_n(s)-\psi_n(t)|_{C(\overline\Omega)}< 1$,
then
  \begin{eqnarray*}
 |\tilde \xi_1|+|\tilde \xi_2| &=& (1-\frac 1{c_\psi}|\psi_n(s)-\psi_n(t)|
   _{C(\overline \Omega)})(|\xi_1|+|\xi_2|) \\
   &\le& (1-\frac 1{c_\psi}
  |\psi_n(s)-\psi_n(t)|_{C(\overline \Omega)})\psi_n(\cdot,s)\\
  &=& \psi_n(\cdot,s)-\frac {\psi_n(\cdot,s)}{c_\psi}
   |\psi_n(s)-\psi_n(t)|_{C(\overline \Omega)} \le \psi_n(\cdot,t).
  \end{eqnarray*}
Hence $\tilde z \in K_n(t)$ and $|\tilde z|_V \le |z|_V$. 
By using extensively this idea we see that (3.1) holds with
$\alpha(t)=c_nt$ for a certain (big) constant $c_n>0$ and $\beta(t)=0$. 
For the detailed proof we refer to [22; Lemma 3.1] or [15; Example 4.5].
\vspace{0.2cm}

(ii) $\{K(t)\} \in \Phi_W$, which is proved by
using the mapping ${\cal F}_\varepsilon(t)z:=
(1-\varepsilon)z$.
In fact, for any $z_n:=[\xi_{1n},\xi_{2n}] \in K_n(t)$ and any small 
$\varepsilon>0$,
we take an integer $N_\varepsilon$ so that 
$|\psi_n -\psi|_{C(\overline Q)}\le \varepsilon c_\psi$ for all 
$n \geq N_\varepsilon$. In this case, we have
 $$(1-\varepsilon)(|\xi_{1n}(x)|+\xi_{2n}(x)|)
  \le (1-\varepsilon)\psi_n(x,t)\le \psi(x,t),
  ~\forall n \geq N_\varepsilon.$$
This shows ${\cal F}_\varepsilon (K_n(t)) \subset K(t)$ for all 
$n \geq N_\varepsilon$. Similarly, ${\cal F}_\varepsilon (K(t)) \subset K_n(t)$ for $n \geq N_\varepsilon$. Hence
$K_n(t) \Longrightarrow K(t)$ on $[0,T]$ and
$\{K(t)\}\in \Phi_W$.\vspace{0.5cm}

As is easily seen from the above examples, the class $\Phi_W$ is strictly 
larger than $\Phi_S$.

Next, we introduce the time-derivative under constraint $\{K(t)\} \in \Phi_W$.
Put
  $$ {\mathcal K}:=\{v \in L^p(0,T;V)~|~ v(t) \in K(t)
  ~{\rm for~a.e.~}t \in [0,T]\} $$
and
 $$ {\mathcal K}_0:=\{ \eta
 \in{\mathcal K} ~|~ \eta'\in 
  L^{p'}(0,T;V^*)\}.$$

\vspace{0.5cm}

\noindent
{\bf Definition 3.2.} Let $\{K(t)\} \in \Phi_W$
and $u_0 \in \overline {K(0)}$.
Then we define an operator $L_{u_0}$ whose graph
$G(L_{u_0})$ is given in
$L^p(0,T;V)\times L^{p'}(0,T;V^*)$ as follows:
$ [u,f] \in G( L_{u_0})$ if and only if
 $$u \in {\cal K},~f \in L^{p'}(0,T; V^*),~~
 \int_0^T \langle \eta'-f, u-\eta\rangle dt 
  \le \frac 12 |u_0-\eta(0)|^2_H, 
  ~~\forall \eta \in 
     {\mathcal K}_0. \eqno{(3.4)}$$
\vspace{0.2cm}

We prove the most important property of 
$L_{u_0}$ in the next theorems.\vspace{0.5cm}

\noindent
{\bf Theorem 3.1.} {\it Let $\{K(t)\} \in \Phi_W$ and $u_0 \in 
\overline{K(0)}$. Then $L_{u_0}$ is maximal monotone from 
$D(L_{u_0})\subset L^p(0,T;V)$
 into 
$L^{p'}(0,T;V^*)$, and the domain $D(L_{u_0})$
is included in the set $\{u\in C([0,T];H)\cap
{\cal K} ~|~u(0)=u_0\}$.}
\vspace{0.5cm}

The characterization and fundamental properties of the mapping $L_{u_0}$
are given in the following theorem.\vspace{0.5cm}

\noindent
{\bf Theorem 3.2.} {\it Let 
$\{K(t)\} \in \Phi_W$. Then we have: 
\begin{description}
\item{(1)} Let $u_0 \in \overline{K(0)}$.
Then $f \in L_{u_0}u$ if and only if there are
$\{\{K_n(t)\}\} \subset \Phi_S$, $\{u_n\} \subset L^p(0,T;V)$ with 
$u_n \in {\cal K}_n:=\{v \in L^p(0,T;V)~|~v(t) \in K_n(t)~for~a.e~t
\in [0,T]\}$ and $u'_n \in 
L^{p'}(0,T;V^*)$, $\{f_n\} \subset 
L^{p'}(0,T;V^*)$ such that 
      $$ K_n(t) \Longrightarrow K(t)~{\rm on~}
      [0,T],\eqno{(3.5)}$$
 $$ u_n \to u~{\rm in~}C([0,T];H)~and~weakly~in~
     L^p(0,T;V), \eqno{(3.6)}$$
 $$ f_n\to f~weakly~in~L^{p'}(0,T;V^*),\eqno{(3.7)} $$
  $$  \int_0^T \langle u'_n-f_n,u_n -v\rangle dt \le 0, ~\forall v \in
      {\mathcal K}_n,~\forall n,\eqno{(3.8)}$$
  $$\limsup_{n\to \infty}\int_{t_1}^{t_2} \langle f_n,u_n\rangle dt
 \le \int_{t_1}^{t_2} \langle f,u \rangle dt, ~\forall t_1,~t_2~with
~0\le t_1\le t_2 \le T.
  \eqno{(3.9)}$$
\item{(2)} 
Let $u_0 \in \overline{K(0)}$
and $f \in L_{u_0}u$. Then, for any $t_1,~t_2 \in 
[0,T]$ with $t_1 \le t_2$,
$$ \int_{t_1}^{t_2}\langle \eta'-f,u-\eta\rangle dt
  + \frac 12|u(t_2)-\eta(t_2)|^2_H \le 
 \frac 12|u(t_1)-\eta(t_2)|^2_H,~~\forall
 \eta \in {\cal K}_0.\eqno{(3.10)}$$
\item{(3)} Let $u_{i0} \in \overline {K(0)},$ and 
$f_i \in L_{u_{i0}}u_i$ for $i=1,2$. Then,
for any $t_1,~t_2 \in [0,T]$ with $t_1\le t_2$, 
 $$ \frac 12 |u_1(t_2)-u_2(t_2)|^2_H \le 
 \frac 12 |u_1(t_1)-u_2(t_2)|^2_H +
  \int_{t_1}^{t_2} \langle f_1-f_2, u_1-u_2 \rangle dt.
  \eqno{(3.11)}$$
\end{description}
}

The proofs of these theorems will be given in the next section.\vspace{1cm}

\noindent
{\bf 4. Proofs of Theorems 3.1 and 3.2}\vspace{0.5cm}

In this section we use the same notation and assume the same assumptions 
as in the previous section.

We introduce another mapping $\tilde L_{u_0}$ whose graph $G(\tilde L_{u_0})$
is given as follows: $f \in \tilde L_{u_0}u$ if and only if $u \in {\cal K}\cap
C([0,T:;H)$ with $u(0)=u_0$, $f \in L^{p'}(0,T;V^*)$ and there exist sequences
$\{\{K_n(t)\}\} \subset \Phi_S$, $\{u_n\} \subset {\cal K}_{n,0}:=\{v \in 
{\cal K}_n~|~ v' \in L^{p'}(0,T;V^*)\}$ and $\{f_n\} \subset L^{p'}(0,T;V^*)$
and (3.5)-(3.9) are fulfilled.

As to the mapping $\tilde L_{u_0}$ we prove:\vspace{0.5cm}

\noindent
{\bf Lemma 4.1.} {\it $\tilde L_{u_0}$ is
a restriction of $L_{u_0}$, namely
$G(\tilde L_{u_0}) \subset G(L_{u_0})$, and
 $$ \int_0^T \langle f-g, u-w \rangle dt \geq 0,
   ~~\forall [u,f] \in G(\tilde L_{u_0}),~~ \forall
   [w,g] \in G(L_{u_0}).
   \eqno{(4.1)}$$
 }
\noindent
{\bf Proof.} Let 
$f \in \tilde L_{u_0}u$ and $\{\{K_n(t)\}\}$, 
$\{u_n\}, ~\{f_n\}$ be sequences as in the definition of 
$f \in \tilde L_{u_0}u$; (3.5)-(3.9) are fulfilled as well. Then, for any $\eta \in {\cal K}_0$, we have
$$ \int_0^T \langle u'_n- f_n, u_n
-{\cal F}_\varepsilon\eta\rangle dt \le 0, $$
since ${\cal F}_\varepsilon \eta \in {\cal K}_{n,0}:=\{v \in {\cal K}_n~|~
v' \in L^{p'}(0,T;V^*)\}$. 
Substituting the expression ${\cal F}_\varepsilon \eta=(1+\varepsilon c_0)
\eta+\varepsilon \sigma_0$ (cf. (3.3)) in the
above inequality and using integration by parts, we get
 $$
  \int_0^T \langle \eta'-f_n, u_n
  -{\cal F}_\varepsilon\eta\rangle dt
+\varepsilon\int_0^T\langle c_0\eta'+\sigma'_0,
  u_n-{\cal F}_\varepsilon \eta \rangle dt 
 \le \frac 12 |u_n(0)-{\cal F}_\varepsilon
   \eta(0)|^2_H
 $$
Now, since ${\cal F}_\varepsilon \eta \to \eta$
in $L^p(0,T;V)\cap C([0,T];H)$ as $\varepsilon \downarrow 0$, we have by letting 
$n\to \infty$ and (3.7) 
 $$ \int_0^T \langle \eta'-f, u-\eta\rangle dt
   \le \frac 12 |u_0-\eta(0)|^2_H,
   ~~\forall \eta\in {\cal K}_0.$$
This implies $f \in L_{u_0}u$. Thus $G(\tilde L_{u_0}) \subset G(L_{u_0})$.

Let $g \in L_{u_0}w$. Then, with the same notation as above,  
it follows from (3.4) that
 $$ \int_0^T \langle({\cal F}_\varepsilon 
 u_n)' -g,
    w -{\cal F}_\varepsilon u_n\rangle dt \le 
  \frac 12|u_0-{\cal F}_\varepsilon(0) u_n(0)|^2_H,
 $$
since ${\cal F}_\varepsilon u_n
=u_n+\varepsilon c_0 u_n+\varepsilon\sigma_0 \in {\cal K}_0$. 
The above inequality is of the form
 $$ \int_0^T \langle u'_n-g +\varepsilon c_0 
  u'_n+\varepsilon \sigma'_0, 
  w-u_n -\varepsilon c_0 u_n-\varepsilon \sigma_0\rangle dt
 ~~~~~~~~~~$$
 $$~~~~~~~ \le \frac 12 |u_0- u_n(0)-\varepsilon 
 c_0u_n(0)-\varepsilon \sigma_0(0)|^2_H.$$  
Hence we derive from this inequality that
 $$
   \int_0^T \langle u'_n-g , 
  w-u_n\rangle dt +\varepsilon\int_0^T\langle 
   c_0u'_n, w \rangle dt ~~~~~~~~~~~~~~~~~~~~~~~~~$$
  $$ \le \frac 12 |u_0- u_n(0)|^2_H+ \varepsilon|(u_0-u_n(0),
 c_0 u_n(0)+\sigma_0(0))_H|+\frac{\varepsilon^2}2|c_0u_n(0)+\sigma_0(0)|^2_H~\eqno{(4.2)} $$
 $$ +\varepsilon  \left | \int_0^T \langle
  c_0 u'_n, u_n+\varepsilon c_0u_n
    +\varepsilon\sigma_0 \rangle dt \right|
 +\varepsilon \left|\int_0^T \langle u'_n- g,
  c_0u_n+\sigma_0 \rangle dt\right |
  ~~~~~~$$
  $$ +\varepsilon \left | \int_0^T \langle\sigma'_0, 
       w-u_n-\varepsilon c_0 u_n
    -\varepsilon \sigma_0\rangle dt \right|.
    ~~~~~~~~~~~~~~~~~~$$
By assumptions (3.6) and (3.7), 
 $$ \sup_{n\geq 1}\left \{
    |f_n|_{L^{p'}(0,T;V^*)} +|u_n|_{L^p(0,T;V)}+|u_n|_{C([0,T];H)}
    \right\} < \infty.$$
Hence, for a constant
$C_1>0$, depending only on $f, u$ and $w$ 
(but independent of $\varepsilon$ and $n$), we infer from (4.2)
that
$$
\int_0^T \langle u'_n-g , 
  w-u_n\rangle dt +\varepsilon c_0\int_0^T\langle 
   u'_n, w \rangle dt
 \le \frac 12 |u_0- u_n(0)|^2_H +\varepsilon
  C_1. \eqno{(4.3)}$$
Also, from (3.8) it follows that
 $$ \int_0^T \langle u'_n - f_n, u_n-
 {\cal F}_\varepsilon w\rangle dt
 =\int_0^T \langle u'_n - f_n, u_n-
   w-\varepsilon c_0 w- \varepsilon
   \sigma_0 \rangle dt \le 0.$$
In a similar way to (4.3) it follows that
 $$\int_0^T \langle u'_n - f_n, u_n-
   w \rangle-\varepsilon c_0 \int_0^T \langle 
    u'_n, w \rangle dt 
   \le \varepsilon C_2, \eqno{(4.4)}$$
with a positive constant $C_2$
depending only on $f, u$ and $w$ (but independent of 
$\varepsilon$ and $n$). Now, adding (4.3)
and (4.4), we arrive 
at an inequality of the form
 $$ \int_0^T \langle f_n -g, u_n -w\rangle dt
    \geq -\varepsilon(C_1+C_2)-\frac 12|u_0-u_n(0)|^2_H,$$
whence (4.1) is obtained by $\varepsilon \downarrow 0$ and 
$n \to \infty$. \hfill $\Box$ 
\vspace{0.5cm}

\noindent
{\bf Corollary 4.1.} {\it $\tilde L_{u_0}$ is monotone from 
$D(\tilde L_{u_0})\subset L^p(0,T;V)$ into $L^{p'}(0,T; V^*)$. Also, if
$u \in D(\tilde L_{u_0})$, then
$u \in {\cal K}\cap C([0,T];H)$ and
$u(0)=u_0$.}\vspace{0.5cm}

Next we show the maximal monotonicity of 
$\tilde L_{u_0}$.\vspace{0.5cm}

\noindent
{\bf Lemma 4.2.} {\it Let $\{K(t)\} \in \Phi_W$
and $u_0 \in \overline {K(0)}$. Then 
$\tilde L_{u_0}$ is
maximal monotone from $L^p(0,T;V)$ into 
$L^{p'}(0,T; V^*)$; more precisely for any 
$f \in L^{p'}(0,T;V^*)$ there exists a unique 
$u\in L^p(0,T;V)$ such that $u(t) \in K(t)$ 
for a.e. 
$t \in [0,T]$ and
  $$ f\in \tilde L_{u_0}u +  Fu.
  \eqno{(4.5)}$$  }
\noindent
\hspace*{-0.2cm}{\bf Proof.} The operator 
$\tilde L_{u_0}+F$ is strictly 
monotone from $L^p(0,T;V)$ into 
$L^{p'}(0,T;V^*)$,
so that the function $u$ satisfying (4.5) is unique.
Now we are to prove the existence of such a function $u$. 

Choose $\{\{K_n(t)\}\} \subset \Phi_S$
such that $K_n(t) \Longrightarrow K(t)$ on
$[0,T]$ and 
$\{f_n\} \subset L^2(0,T;H)$ 
such that $f_n \to f$ in $L^{p'}(0,T;V^*)$ with $|f_n|_{L^{p'}(0,T;V^*)}
\le |f|_{L^{p'}(0,T;V^*)}+1$. 
Also, it is easy to take a sequence 
$\{u_{0,n}\}$ in $V$ such that 
$u_{0,n} \in K_n(0)$ and $u_{0,n} \to u_0$ in 
$H$. For these data we consider
the Cauchy problems
 $$ u'_n(t)+U_n(t)=f_n(t),~U_n(t) \in
   \partial \varphi^t_{K_n}(u_n(t))~{\rm in~}H,
    ~{\rm a.e.~}t \in [0,T],
   u_n(0)=u_{0,n}, \eqno{(4.6)}$$
where $\varphi^t_{K_n}(\cdot)$ is a 
time-dependent proper, l.s.c. and convex
functions on $H$ given by (3.2) with $K(t)$ replaced by $K_n(t)$. 
As was mentioned in section 3, (4.6) possesses a
unique solution $u_n \in W^{1,2}(0,T;H)$ such that $t \to \varphi^t_{K_n}(u_n(t))$ is absolutely continuous on $[0,T]$ 
(hence, $u_n \in {\mathcal K}_{n,0}$ and $u_n(t) \in K_n(t)$ for all 
$ t \in [0,T]$).

We note here that 
 $$ U_n(t) \in \partial_*I_{K_n(t)}(u_n(t))+Fu_n(t),~~{\rm a.e.~}t \in [0,T],
  \eqno{(4.7)}$$
since $\partial \varphi^t_{K_n}(z) \subset \partial_* \varphi^t_{K_n}(z) 
= \partial_*I_{K_n(t)}(z)+Fz$ for all $z \in K_n(t)$ (cf. [10; Theorem 2] or
[20; Theorem 5.2]), and by the expression
 $$
 {\cal F}^2_\varepsilon(t) z:={\cal F}_\varepsilon(t)({\cal F}_\varepsilon(t)z)
    = z +2\varepsilon c_0 z+
       \varepsilon^2 c_0^2 z
       +\varepsilon^2 c_0 \sigma_0(t)
     +2\varepsilon \sigma_0(t) \eqno{(4.8)} $$
that ${\cal F}^2_\varepsilon(t)K_m(t)
\subset K_n(t)$ for any large $n,~m$ and small $\varepsilon >0$. 
Therefore, from (4.6)-(4.8) it
 follows for any $m,~n$ and $t \in [0,T]$ that
 $$\int_0^t(u_n', u_n-{\cal F}^2_\varepsilon u_m)_H d\tau
 +\int_0^t\langle  Fu_n,u_n-{\cal F}^2_\varepsilon u_m \rangle d\tau \le
   \int_0^t \langle f_n,u_n-{\cal F}^2_\varepsilon u_m\rangle d\tau,
   \eqno{(4.9)}$$
and similarly
$$\int_0^t(u_m', u_m-{\cal F}^2_\varepsilon u_n)d\tau
 +\int_0^t\langle  Fu_m,u_m-{\cal F}^2_\varepsilon u_n \rangle 
   d\tau\le 
   \int_0^t \langle f_m,u_m-{\cal F}^2_\varepsilon 
    u_n\rangle d\tau.
    \eqno{(4.10)}$$
We observe from the energy estimate for (4.6)
 that
$\{u_m\}$ is bounded in $L^p(0,T;V)$ and
$C([0,T];H)$, namely there is a constant
$C_3>0$ such that
 $$ |u_m|_{L^p(0,T;V)}+ |u_m|_{C([0,T];H)}
    \le C_3,~~\forall m.\eqno{(4.11)}$$
In fact, fixing a number $n$, we get an estimate
of the form (4.11) by applying the Gronwall's
inequality to (4.10). 

Now substitute (4.8) with $z=u_m$ and $z=u_n$
in (4.9) and (4.10), respectively. Then the sum of resultants 
gives an inequality of the form
 $$ \int_0^t(u'_n-u'_m, u_n-u_m)_H d\tau
 + \int_0^t \langle Fu_n-Fu_m,u_n-u_m\rangle d\tau $$
 $$ ~~~~~~~\le \int_0^t \langle f_n-f_m,u_n-u_m
 \rangle d\tau +\varepsilon (D_{1,n,m}(t)+D_{2,n,m}(t)),$$ 
where 
 $$
D_{1,n,m}(t):=(2+ \varepsilon c_0) 
\int_0^t \{\langle u'_n, c_0u_m+\sigma_0\rangle
 +\langle u'_m, c_0 u_n+\sigma_0\rangle\} d\tau 
 $$
and 
 $$
  D_{2,n,m}(t):=(2+\varepsilon c_0) \int_0^t \{\langle  Fu_n-f_n,c_0 u_m
   +\sigma_0 \rangle 
   + \langle  Fu_m-f_m,c_0 u_n+\sigma_0 
   \rangle\} d\tau. 
  $$
Here, by integration by parts, $D_{1,n,m}$ is
arranged in the form: 
 \begin{eqnarray*}
 D_{1,n,m}(t) &=& (2+\varepsilon c_0)c_0
 \{(u_n(t),u_m(t))-(u_n(0),u_m(0))\} \\
&&+ (2+\varepsilon c_0)\{(u_n(t)+u_m(t),\sigma_0(t))-(u_n(0)+u_m(0),\sigma_0(0))\}
\\
  && -(2+\varepsilon c_0) \int_0^t \langle \sigma'_0, u_n+u_m \rangle d\tau
 \end{eqnarray*}
for all $t \in [0,T]$, so that  
it follows from (4.11) that $D_{1,n,m}(t)$ is dominated by a 
positive constant
$D_1$ independent of $m,n$ and $t \in [0,T]$; for instance, 
$D_1:=6c_0C_3^2+12C_3
\{|\sigma_0|_{C([0,T];V)}+|\sigma'_0|_{L^{p'}(0,T;V^*)}\}$ for all small $\varepsilon>0$.
Similarly, it follows from (4.11) that 
 $$ D_{2,n,m}(t) \le D_2:=12(C^p_3+ |f|_{L^{p'}(0,T;V^*)}+1)
   (c_0C_3+|\sigma_0|_{L^p(0,T;V)}).$$
Therefore, putting $C_4:=D_1+D_2$, we obtain
 $$ \frac 12 |u_n(t)-u_m(t)|^2_H
   + \int_0^t\langle Fu_n-Fu_m,u_n-u_m\rangle d\tau$$
  $$ \le  \frac 12 |u_{0,n}-u_{0,m}|^2_H
    +\int_0^t\langle f_n-f_m,u_n-u_m
 \rangle d\tau+\varepsilon C_4.$$
By letting $n,m \to \infty$ and $\varepsilon \downarrow 0$ in the above
inequality we see
that for some $u \in C([0,T];H)\cap L^p(0,T;V)$
such that $ u_n \to u~{\rm in~}C([0,T];H)~ {\rm and~weakly~in~} L^p(0,T;V)$
and
 $$  0 \le \int_{t_1}^{t_2} \langle Fu_n-Fu_m,u_n-u_m\rangle dt \le
 \int_0^T \langle Fu_n-Fu_m,u_n-u_m\rangle dt
 \to 0~~{\rm as~}n,m\to, \infty \eqno{(4.12)}$$
for all $0\le t_1 \le t_2 \le T$. Taking a subsequence of $\{n\}$
if necessary, we have 
  $$ Fu_n \to \ell^*~{\rm weakly~in~}L^{p'}(0,T;V^*)~{\rm for~some~}
     \ell^* \in L^{p'}(0,T; V^*),$$ 
which imply by (4.12)
 $$  \limsup_{n \to \infty} \int_{t_1}^{t_2}
  \langle  Fu_n, u_n \rangle dt \le 
   \int_{t_1}^{t_2}\langle \ell^*, u \rangle 
   d t,~ {\rm for~all~} 0\le t_1\le t_2\le T.$$
From the maximal monotonicity of $F$ we obtain that $\ell^*=Fu$ and
$$ \lim_{n\to \infty}\int_{t_1}^{t_2} \langle Fu_n,u_n\rangle dt=\int_{t_1}^{t_2} \langle 
Fu,u \rangle dt,~{\rm for~all~} 0\le t_1\le t_2\le T.
 $$ 
Consequently,
$ f_n - Fu_n \to f- Fu~{\rm weakly~in~} L^{p'}(0,T;V^*)$ and   
 $$\lim_{n\to \infty}\int_{t_1}^{t_2} \langle f_n-  Fu_n,u_n\rangle dt = \int_{t_1}^{t_2} 
\langle f-  Fu,u\rangle dt, ~{\rm for~all~}
  0\le t_1\le t_2\le T.$$ 
Therefore, by definition,
$f-  Fu \in \tilde L_{u_0}u$. 

Thus we have seen that the range of $\tilde L_{u_0}+ F$ is the whole of 
$L^{p'}(0,T;V^*)$. Since
$\tilde L_{u_0}$ is monotone by Corollary 4.1, we conclude that it is maximal monotone from $L^p(0,T;V)$ into $L^{p'}(0,T;V^*)$.
From the definition of 
$\tilde L_{u_0}$ we see that $D(\tilde L_{u_0})
\subset \{u \in C([0,T];H)\cap {\cal K},~
u(0)=u_0\}$. \hfill $\Box$ 
\vspace{0.5cm}

Now it follows from the Lemma 4.2 that $\tilde L_{u_0}= L_{u_0}$, and 
Theorem 3.1 is obtained. 
\vspace{0.5cm}

\noindent
{\bf Proof of Theorem 3.2.} (1) is obtained from the fact that
$\tilde L_{u_0}=L_{u_0}$. Next, we prove (2).
Corresponding to $f \in L_{u_0}u$, choose 
$\{\{K_n(t)\}\} \subset \Phi_S$, $\{u_n\}$ with $u_n \in {\cal K}_{n,0}$ and 
$\{f_n\} \subset L^{p'}(0,T;V^*)$
so that conditions (3.5)-(3.9) hold. Take as a test function $v$ in (3.8)
  $$v:=\left \{ 
        \begin{array}{l}
          {\cal F}_\varepsilon \eta~~~{\rm on~}[t_1,t_2],\\
          u_n~~~~~{\rm on~}[0,t_1)\cup (t_2,T],
        \end{array} \right. $$
for any $ \eta \in {\cal K}_0$ and small $\varepsilon>0$ to obtain
  $$\int_{t_1}^{t_2} \langle u'_n-f_n, u_n-
  {\cal F}_\varepsilon \eta \rangle dt\le 0.
  $$
Applying integration by parts to this
inequality and substitute the expression of
${\cal F}_\varepsilon \eta$ in it, we see that
 $$ \int_{t_1}^{t_2} \langle \eta'
 +\varepsilon c_0\eta'
 +\varepsilon\sigma'_0
 -f_n, u_n- {\cal F}_\varepsilon \eta \rangle dt
 +\frac 12 |u_n(t_2)-{\cal F}_\varepsilon(t_2)\eta(t_2)|^2_H 
\le \frac 12 |u_n(t_1)-{\cal F}_\varepsilon(t_1)\eta(t_1)|^2_H.$$
Therefore, letting $n\to \infty$ and 
$\varepsilon\to 0$ yield (3.10).

Next we show (3.11). Choose
$\{K_{1,n}\}\in \Phi_S$, $\{u_{1,n}\}$ and $\{f_{1,n}\}$ so that 
conditions (3.5)-(3.9) hold corresponding to $f_1 \in L_{u_{10}}u_1$ as well as
$\{K_{2,m}\}\in \Phi_S$, $\{u_{2,m}\}$ and $\{f_{2,m}\}$
corresponding to $f_2 \in L_{u_{20}}u_2$.
Noting that ${\cal F}^2_\varepsilon u_{1,n}
\in {\cal K}_{2,m}$ as well as ${\cal F}^2_\varepsilon u_{2,m}
\in {\cal K}_{1,n}$ for all large $n,~m$, we observe by taking
    $$ v:=\left \{
         \begin{array}{l}
          {\cal F}^2_\varepsilon u_{1,n}~~~~~{\rm on~}[t_1,t_2],\\
          u_{2,m}~~~~~~~~{\rm on~}[0,t_1)\cup (t_2,T],
         \end{array} \right. $$
in (3.8) for $u_{2,m}$ that
 $$ \int_{t_1}^{t_2} \langle u'_{2,m}-f_{2,m}, u_{2,m}
 -{\cal F}^2_\varepsilon u_{1,n} \rangle dt
  \le 0. \eqno{(4.13)}$$
Similarly, for large $n,~m$,
 $$ \int_{t_1}^{t_2} \langle u'_{1,n}-f_{1,n}, u_{1,n}
 -{\cal F}^2_\varepsilon u_{2,m} \rangle dt \le 0.
 \eqno{(4.14)}$$
Substitute the expression of 
${\cal F}^2_\varepsilon u_{1,n}$ and ${\cal F}^2_\varepsilon u_{2,m}$
in (4.13) and (4.14) and add them to get
$$\int_{t_1}^{t_2} \langle u'_{1,n}-u'_{2,m},
   u_{1,n}-u_{2,m}\rangle dt
   \le  \int_{t_1}^{t_2} \langle f_{1,n}-f_{2,m},u_{1,n}-u_{2,m}\rangle dt
   + \varepsilon C_5,$$
where $C_5$ is a positive constant independent $t_1,~t_2,~\varepsilon$ and $n,~m$.
Hence,
$$\frac 12 |u_{1,n}(t_2)-u_{2,m}(t_2)|^2_H\le 
   \frac 12 |u_{1,n}(t_1)-u_{2,m}(t_1)|^2_H 
   +\int_{t_1}^{t_2} \langle f_{1,n}-f_{2,m},
u_{1,n}-u_{2,m}\rangle dt+ \varepsilon C_5,$$
and we obtain (3.11) by passing to the limit
$n,~m\to \infty$ and $\varepsilon \to 0$. \hfill $\Box$ \vspace{0.5cm}

\noindent
{\bf Remark 4.1.} In Hilbert spaces similar operators to
$L_{u_0}$ were considered in the
time-independent case $K=K(t)$ (cf. [8]) and 
it was generalized to the time-dependent case $K(t)$ (cf. [21]). 
In the Banach space set-up (cf. [19]), the similar results were discussed,
too.

\noindent
{\bf Remark 4.2.} Theorem 3.1 gives a
generalization of the results of [19, 21] in a class of weak variational inequalities. Moreover
it is expected to compose $L_{u_0}$ for various  constraint set $K(t)$ in a much wider class 
$\Phi_W$ than in this paper, for instance the
classes specified in [15, 23].
\vspace{1cm}

\noindent
{\bf 5. Perturbations of semimonotone type}\vspace{0.5cm}

We assume that $H,~V$ and $W$ be the same as in section 2; $V$ is 
dense in $H$ with compact injection and $W$ is separable and dense in $V$
with continuous injection, and $1<p<\infty, \frac 1p+\frac 1{p'}=1$. 

Let $A(t,v,u)$ be a singlevalued mapping from 
$[0,T]\times H \times V$ into $V^*$, and assume that: 
\begin{description}
\item{(a)} (Boundedness) There are positive constants 
$c_1,~c_2$ such that
 $$ |A(t,v,u)|_{V^*} \le c_1|u|^{p-1}_V +c_2,~
   ~\forall (t, v, u) \in [0,T] \times H \times V.$$
\item{(b)} (Coerciveness) There are positive
constants $c_3,~c_4$ such that
  $$\langle A(t,v,u), u \rangle \geq c_3|u|^p_V-c_4,
        ~~\forall (t, v, u) \in [0,T] \times H \times V.$$
\item{(c)} (Semimonotonicity) For each 
$v \in H$ and $t \in [0,T]$, 
the mapping $u \to A(t,v,u)$ is demicontinuous
from $D(A(t,v,\cdot))=V$ into $V^*$ and monotone, namely 
   $$ \langle A(t,v,u_1)-A(t,v,u_2), u_1- u_2  \rangle \geq 0,
      ~~\forall u_1,~u_2 \in V, $$
Moreover, for each $u \in V$
the mapping $(t,v) \to A(t,v,u)$ is continuous from $[0,T]\times H$
into $V^*$.
\end{description}

We derive some properties of $A(t,v,u)$ from the
above conditions. 
\vspace{0.5cm}

\noindent
{\bf Lemma 5.1.} {\it Assume that $t_n \to t$ in $[0,T]$, $v_n \to v$
in $H$, $u_n \to u$ weakly in $V$ 
and
   $$ A(t_n,v_n,u_n) \to \alpha^*~ {\it weakly~in~} V^*,
       ~~\limsup_{n\to \infty} \langle A(t_n,v_n,u_n),u_n\rangle
         \le \langle \alpha^*, u\rangle. \eqno{(5.1)}$$
Then $\alpha^* =A(t,v,u)$ and $\lim_{n\to \infty} \langle 
A(t_n,v_n,u_n),u_n\rangle
   = \langle \alpha^*, u\rangle$.  }\vspace{0.3cm}

\noindent
{\bf Proof.} By condition (c), for each $t \in [0,T]$ and $v\in H$, the
mapping $u \to A(t,v,u)$ is monotone and
demicontinuous from $V$ into $V^*$ and hence
maximal monotone (cf. [20; Theorem 4.2]. Let 
$\eta$ be any
element in $V$. Then, by assumption (5.1) with (c),
 $$
  \langle \alpha^*-A(t,v,\eta), u-\eta \rangle
 \geq limsup_{n\to \infty} \langle A(t_n, v_n,u_n)- A(t_n,v_n,\eta),
     u_n-\eta\rangle \geq 0. $$
The maximal monotonicity of $\eta \to 
A(t,v,\eta)$ implies that 
$\alpha^*=A(t,v,u)$ as well as $\langle
A(t_n,v_n,u_n), u_n\rangle \to \langle \alpha^*, u \rangle$.
\hfill $\Box$\vspace{0.5cm}

The above lemma ensures that $A(t,v,u)$ is continuous from $[0,T]\times H\times
V$ into $V^*_w$, so that for every $v \in L^p(0,T;H)$ and 
$u\in L^p(0,T;V)$ the function
$t \to A(t,v(t),u(t))$ is weakly measurable on $[0,T]$ and hence strongly 
measurable on $[0,T]$ thanks to the separability of $V^*$. As a consequence,
by condition (a) we have
 $A(\cdot,v,u) \in L^{p'}(0,T; V^*)$ for
every $v \in L^p(0,T;H)$ and $u\in L^p(0,T;V)$. \vspace{0.5cm}

Now we introduce the mapping ${\cal A}: D({\cal A})= L^p(0,T; H)\times 
L^p(0,T;V) \to L^{p'}(0,T;V^*)$ by putting
  $$ [{\cal A}(v,u)](t):= A(t,v(t),u(t))~~{\rm for~a.e.~}t \in [0,T].
\eqno{(5.2)}$$

\noindent
{\bf Lemma 5.2.} {\it ${\cal A}$ satisfies the following properties:
\begin{description}
\item{(i)} (Boundedness) $|{\cal A}(v,u)|_{L^{p'}(0,T;V^*)} 
   \le c_1|u|_{L^p(0,T;V)}^{\frac p{p'}} +c_2T^{\frac 1{p'}} $
for all $v \in L^p(0,T;H)$ and $u\in L^p(0,T;V)$.
\item{(ii)} (Coerciveness) $\int_0^T \langle {\cal A}(v,u), u\rangle \geq
c_3 |u|_{L^p(0,T;V)}^p - c_4T$ 
for all $v \in L^p(0,T;H)$ and $u\in L^p(0,T;V)$.
\item{(iii)} (Semimonotonicity) ${\cal A}$ is demicontinuous from 
$L^p(0,T;H)\times L^p(0,T;V)$ into $L^{p'}(0,T;$ $V^*)$ and for every 
$v \in L^p(0,T;H)$ the mapping
$u\to {\cal A}(v,u)$ is monotone from $L^p(0,T;V)$ into $L^{p'}(0,T;V^*)$,
namely,
 $$ \int_0^T \langle {\cal A}(v, u_1) -{\cal A}(v,u_2), u_1-u_2 \rangle dt
    \geq 0,~~\forall u_1,~u_2 \in L^p(0,T;V).$$
\item {(iv)} (Continuity in $v$) For every $u \in L^p(0,T;V)$, the mapping 
$v \to {\cal A}(v,u)$
is continuous from $L^p(0,T;H)$ into $L^{p'}(0,T;V^*)$.
\item{(v)} Assume that $v_n \to v$ in $L^p(0,T;H)$, $u_n \to u$
weakly in $L^p(0,T;V)$, ${\cal A}(v_n,u_n) \to \alpha^*$ in $L^{p'}(0,T;V^*)$
and 
  $$ \limsup_{n\to \infty}\int_0^T\langle {\cal A}(v_n,u_n),u_n\rangle dt
      \le \int_0^T\langle \alpha^*, u\rangle dt.  $$ 
Then, 
   $$\alpha^*={\cal A}(v,u),~~~\lim_{n\to \infty} \int_0^T\langle 
   {\cal A}(v_n,u_n), u_n \rangle dt=\int_0^T\langle 
   {\cal A}(v,u), u \rangle dt. $$ 
\end{description}
}

The statements (i)-(iv) of Lemma 5.2 are straightforwardly obtained from 
conditions (a), (b) and (c), and (v) is proved in the same way just as 
 Lemma 5.1.
 \vspace{0.5cm}

We are now in a position to state a perturbation result of $L_{u_0}$.
\vspace{0.5cm}

\noindent
\noindent
{\bf Theorem 5.1.} {\it Let ${\cal A}:={\cal A}(v,u)$ be given by (5.2).
Let $\{K(t)\} \in \Phi_W$ and $u_0 \in \overline{K(0)}$ and assume that
there is a positive number $\kappa >0$ such that
  $$ \kappa B_W(0) \subset K(t),~~\forall t \in [0,T]. \eqno{(5.3)}$$
Then, for any 
$f  \in L^{p'}(0,T;V^*)$ there exists a function $u \in D(L_{u_0})$ such
that
    $$ f \in L_{u_0}u+{\cal A}(u,u) \eqno{(5.4)}$$
 }
\noindent
{\bf Proof.}
By Lemma 5.2, for each
$v \in L^p(0,T;V)$ the mapping 
$u \to {\cal A}(v,u)$ is
maximal monotone from $L^p(0,T;V)$ into $L^{p'}(0,T;V^*)$, bounded and 
coercive.
Therefore, on account of the well-known result of maximal monotone 
perturbations (cf. [4, 10, 11, 20]) the range of $L_{u_0}+{\cal A}(v,\cdot)$ 
is the whole of
$L^{p'}(0,T;V^*)$, namely for any $f \in L^{p'}(0,T;V^*)$ there is an
element $u \in D(L_{u_0})$ satisfying 
 $$ f \in L_{u_0}u+{\cal A}(v,u), \eqno{(5.5)}$$
and $u$ is unique by (3) of Theorem 3.2. Now, we denote by 
${\cal S}$ the mapping which
assigns to each $v \in L^p(0,T;V)$ the unique solution $u \in L^p(0,T;V)$
of (5.5), namely $u={\cal S}v$. On account of (2) of
Theorem 3.2,
 $$ \int_0^T \langle \eta' -(f-{\cal A}(v,u)),
 u-\eta \rangle dt \le \frac 12 
 |u_0-\eta(0)|^2_H,~~\forall \eta \in 
 {\cal K}_0, \eqno{(5.6)}$$
where ${\cal K}_0=\{\eta\in L^p(0,T;V)~|~
\eta' \in L^{p'}(0,T;V^*),~\eta(t) \in K(t)~
{\rm for~a.e.~}t \in [0,T]\}$.

Next, we show the uniform estimate for solutions $u$ of (5.5). To do so,
use (3) of Theorem 3.2 for $0 \in L_00$ and $f-{\cal A}(v,u) \in L_{u_0}u$.
Then we get
 $$ \frac 12 |u(t)|^2_H 
 +\int_0^t \langle {\cal A}(v,u),u\rangle d\tau
  \le \frac 12 |u_0|^2_H+\int_0^t \langle f, u\rangle d\tau,~~\forall t\in [0,T]. $$
By condition (b), this inequality implies that
 $$ |u(t)|^2_H 
     +2c_3 \int_0^t |u|^p_V d\tau- 2c_4t \le 
     |u_0|^2_H
     +2 \int_0^t |f|_{V^*}|u|_Vd\tau ,
     ~~\forall t 
     \in [0,T], $$
and by Young's inequality
  $$ |u(t)|^2_H +c_3 \int_0^t |u|^p_Vd\tau 
  \le 2c_4T+|u_0|^2_H
   + c_5\int_0^T|f|^{p'}_{V^*}d\tau,
   ~~\forall t \in [0,T],$$
where $c_5$ is a positive constant, for instance
$ c_5:= \frac{2^{p'}}{p'c_3^{p'-1}p^{p'-1}}$. 
We derive from the above estimate that 
  $$ \sup_{0 \le t \le T}|u(t)|_H + |u|_{L^p(0,T;V)} \le M_2, \eqno{(5.7)}$$
$M_2$ being a positive constant, for instance
$M_2:=2\left\{1+ 2c_4T+|u_0|^2_H
   + c_5\int_0^T|f|^{p'}_{V^*}dt \right\}$. 
Moreover, with the same notation as above, we have by condition (a)
   \begin{eqnarray*}
  |f-{\cal A}(v,u)|_{L^{p'}(0,T;V^*)}
   &\le&  |f|_{L^{p'}(0,T;V^*)}
      +|{\cal A}(v,u)|_{L^{p'}(0,T;V^*)} \\
   &\le&  |f|_{L^{p'}(0,T;V^*)}
    +c_1|u|^{\frac p{p'}}
    _{L^p(0,T;V)}+c_2T^{\frac 1{p'}}=:M_3,
   \end{eqnarray*}
namely it follows that for all 
$v \in L^p(0,T;V)$ and $u ={\cal S}v$
 $$ \inf_{\ell^* \in L_{u_0}u}
    |\ell^*|_{L^{p'}(0,T;V^*)} \le M_3
          \eqno{(5.8)}$$
We take as a constant $M_0$ of Theorem 2.1 the sum $M_2+(1+M_2+c_6)M_3$, where
$c_6$ is a positive constant satisfying
$|\cdot|_{L^1(0,T;W^*)}\le c_6 |\cdot|_{L^{p'}
(0,T;V^*)}$, and
consider the set 
 $$ {\cal X}_0:={\rm conv}[Z(\kappa,M_0,u_0)].$$
By Theorem 2.1, it is non-empty, closed, convex
and compact in $L^p(0,T;H)$.
We are going to apply the Schauder fixed point theorem to ${\cal S}$
in ${\cal X}_0$.

First we check that ${\cal S}({\cal X}_0) \subset {\cal X}_0$. In fact,
for each $v \in {\cal X}_0$ the solution $u:={\cal S}v$ of (5.5) satisfies (5.6) and
estimate (5.7)-(5.8), so that 
$u \in Z(\kappa, M_0, u_0) \subset {\cal X}_0$. Thus ${\cal S}$ maps ${\cal X}_0$ into itself.
Next we show the continuity of ${\cal S}$ in ${\cal X}_0$ with respect to the
topology of $L^p(0,T;H)$. Assume that
$v_n \in {\cal X}_0$ and $v_n \to v$ in $L^p(0,T;H)$ 
(as $n \to \infty$). Clearly, $v_n \to v$ weakly in $L^p(0,T;V)$
and $v \in {\cal X}_0$. 
Putting $u_n={\cal S}v_n$, we 
infer from the compactness of ${\cal X}_0$ that $\{v_{n_k}\}$ and $\{u_{n_k}\}$ converge
in $L^p(0,T;H)$ and in $H$ a.e. on $[0,T]$
to some functions $v$ and $u \in {\cal X}_0$, respectively, for a certain 
subsequence $\{n_k\}$ of $\{n\}$. In this case, with the notations
  $$ f = \ell^*_{n_k} + {\cal A}(v_{n_k}, u_{n_k}),
  ~ \ell^*_{n_k} \in L_{u_0}u_{n_k}, $$
we may assume by the boundedness of ${\cal A}$ that
 $$ {\cal A}(v_{n_k},u_{n_k})  \to \alpha^* ~{\rm weakly~in~}L^{p'}
(0,T;V^*) ~{\rm for~ some~}\alpha^* \in L^{p'}(0,T;V^*),$$
and 
 $$ \lim_{k\to \infty}\int_0^T \langle {\cal A}(v_{n_k},u_{n_k}), 
u_{n_k}-u \rangle dt  ~~{\rm exists.}$$
In order to prove that $\ell^*:=f -\alpha^* \in L_{u_0}u$, we observe
 \begin{eqnarray*}
    0&=& \lim_{k \to \infty}\int_0^T \langle f, u_{n_k}-u \rangle dt \\
  &=&\lim_{k \to \infty}\left\{\int_0^T \langle \ell^*_{n_k}, u_{n_k}
     -u \rangle dt +\int_0^T\langle {\cal A}(v_{n_k},u_{n_k}), u_{n_k}-u \rangle dt
       \right\}\\
  &=&\lim_{k \to \infty}\int_0^T \langle \ell^*_{n_k}, u_{n_k}-u \rangle dt
    +\lim_{k \to \infty}\int_0^T \langle {\cal A}(v_{n_k},u_{n_k}), 
     u_{n_k}-u \rangle dt.
 \end{eqnarray*}
Now we show $\lim_{k \to \infty}\int_0^T \langle \ell^*_{n_k}, u_{n_k}-u 
\rangle dt \le 0$ by contradiction.
Assuming that 
 $$\lim_{k \to \infty}\int_0^T \langle \ell^*_{n_k}, u_{n_k}-u \rangle dt >0,$$
we have by (iii) and (vi) of Lemma 5.2
  \begin{eqnarray*}
 0 &>& \lim_{k \to \infty}\int_0^T \langle {\cal A}(v_{n_k},u_{n_k}), 
       u_{n_k}-u \rangle dt\\   
   &=& \lim_{k \to \infty}\left\{\int_0^T \langle {\cal A}(v_{n_k},u_{n_k})
     -{\cal A}(v_{n_k},u), u_{n_k}-u \rangle dt
 + \int_0^T\langle {\cal A}(v_{n_k},u), u_{n_k}-u \rangle dt \right\} \\
   &\geq& \limsup_{k\to \infty}\int_0^T \langle 
    {\cal A}(v_{n_k},u), 
       u_{n_k}-u \rangle dt =0, 
  \end{eqnarray*}
which is a contradiction. Consequently, we have
 $$\lim_{k\to\infty} \int_0^T \langle \ell^*_{n_k},u_{n_k}-u\rangle dt \le 0.$$
as well as
$ \ell^*_{n_k}=f-{\cal A}(v_{n_k},u_{n_k}) \in L_{u_0}u_{n_k},~u_{n_k} \to u
      ~{\rm weakly~in~}L^p(0,T;V)$ and 
$\ell^*_{n_k} \to \ell^*:=f- \alpha^*~{\rm weakly~in~}L^{p'}(0,T;V^*).$
By the maximal monotonicity of $L_{u_0}$, we obtain that
$\ell^*=f-\alpha^* \in L_{u_0}u$, namely $f \in L_{u_0}u +\alpha^*$, and
$\lim_{k\to\infty} \int_0^T \langle \ell^*_{n_k},u_{n_k}\rangle dt=
   \int_0^T \langle \ell^*,u \rangle dt$, or equivalently 
$\lim_{k\to\infty} \int_0^T \langle {\cal A}(v_{n_k},u_{n_k}),u_{n_k}
     \rangle dt =
   \int_0^T \langle \alpha^*,u \rangle dt$. Therefore it is concluded by
(v) of Lemma 5.2 that $\alpha^*={\cal A}(v,u)$ and $u={\cal S}v$.
By the uniqueness of solution to (5.5), we see that $u_n \to u$ in $L^p(0,T;H)$
without extracting any subsequence from $\{u_n\}$.
Thus ${\cal S}$ is continuous
in ${\cal X}_0$ in the topology of $L^p(0,T;H)$, so that ${\cal S}$ possesses
at least one fixed point $u:={\cal S}u$, which gives a solution of (5.4).
\hfill $\Box$ \vspace{0.5cm}

\noindent
{\bf (Application 1)} \vspace{0.3cm}

The model problem mentioned in the introduction
is here discussed precisely in our framework.

Let $\Omega$ be a bounded and smooth domain in ${\bf R}^N$ and
$Q:=\Omega\times (0,T),~0<T<\infty$.
We use our abstract theorems in the set-up 
 $$H:=L^2(\Omega)\times L^2(\Omega),~
    V:=H^1_0(\Omega)\times H^1_0(\Omega),~
    W:=W_0^{1,q}\times W_0^{1,q},~\max\{N, 2\}<q<\infty.$$
Hence $V^*=H^{-1}(\Omega)\times 
H^{-1}(\Omega)$, $V\subset H\subset V^*
\subset W^*$ and 
$W \subset C(\overline{\Omega})\times C(\overline{\Omega})$ with compact 
embeddings.

Let $\psi=\psi(x,t)$ be an obstacle function prescribed in $C(\overline Q)$ 
so that $\psi \geq c_\psi$ on $\overline Q$ for a positive constant 
$c_\psi$, and define a constraint set $K(t)$ by
 $$ K(t):=\{[\xi_1,\xi_2] \in V~|~
 |\xi_1|+|\xi_2| \le \psi(\cdot,t)~{\rm a.e.~in~}
 \Omega\},~~\forall t\in [0,T].$$
Then, by virtue of Theorem 3.1, the maximal monotone mapping
$L_{u_0}: D(L_{u_0}) \subset L^2(0,T;V)\to
 L^2(0,T;V^*)$ is well defined for any given initial datum
$u_0:=[u_{10},u_{20}] \in \overline {K(0)}$.

Now, we define a nonlinear mapping 
$A(t,v,u): [0,T]\times H\times V \to V^*$ by
 $$ \langle A(t,v,u), \xi\rangle :=
 \int_\Omega \{a_1(x,t,v)\nabla u_1\cdot\nabla \xi_1
 +a_2(x,t,v)\nabla u_2\cdot \nabla \xi_2\}dx,
     $$
  $${\rm for~} v:=[v_1,v_2]\in H, 
  ~u:=[u_1,u_2] \in V,~
     \xi=[\xi_1,\xi_2] \in V,~t \in [0,T], $$
where $a_1(x,t,v)$ and $a_2(x,t,v)$ are functions satisfying the
Carath\'eodory condition on $\overline \Omega \times [0,T]\times {\bf R}^2$
and
 $$ a_*\le a_i(x,t,v) \le a^*~~{\rm for~a.e.~} (x,t) \in \Omega 
 \times (0,T)~{\rm and ~all~}v \in {\bf R}^2,~ i=1,2,$$
for positive constants $a_*,~a^*$. Under the
above assumptions, we easily check the
conditions (a), (b) and (c) in section 5 as well as condition (5.3)
of Theorem 5.1 by the strict positiveness of 
$c_\psi$.
Accordingly we can apply Theorem 5.1
to solve our model problem for given data
$u_0:=[u_{01},u_{02}] \in \overline{K(0)}$ 
and $f =[f_1,f_2]\in L^2(0,T; V^*)$ in the form 
$ f \in L_{u_0}u+{\cal A}(u,u)$.
This functional inclusion is equivalent to
the following weak variational form: 
 $$ u:=[u_1,u_2] \in L^2(0,T;V)\cap C([0,T]; H),
    |u_1|+|u_2|\le \psi~{\rm a.e.~on~}Q,
     ~u(0)=u_0; $$ 
 $$ \int_Q \{\xi_{1,t}(u_1-\xi_1)
    + \xi_{2,t} (u_2-\xi_2)\}dxdt ~~~~ ~~~~~~~~~~~~~~~~~~~~~~~~~~~~~~~~~~~~~~~~~~~~$$ 
  $$ +\int_Q \{a_1(x,t,u)\nabla u_1\cdot 
\nabla (u_1-\xi_1)
    +a_2(x,t,u)\nabla u_2\cdot \nabla (u_2-
  \xi_2)\}dxdt $$
  $$ \le \int_0^T\langle f, u-\xi \rangle dt 
  +\frac 12 \{|u_{10}-\xi_1(0)|^2_{L^2(\Omega)}
  +|u_{20}-\xi_2(0)|^2_{L^2(\Omega)}\}, $$
  $$ \forall \xi=[\xi_1,\xi_2] \in L^2(0,T;V)
   \cap W^{1,2}(0,T; H),~|\eta_1|+|\eta_2|
   \le \psi~ {\rm~a.e.~on~}Q.$$
\vspace{0.3cm}

\noindent
{\bf (Application 2)}\vspace{0.3cm}

Finally, we consider parabolic variational
inequalities with gradient constraints.

Let $\Omega$ be a bounded and smooth domain in ${\bf R}^N$
and put
 $$ H:=L^2(\Omega),~V:=H^1_0(\Omega),~
    W:=W^{2,q}_0(\Omega), ~\max\{N,2\} <q <\infty.$$
For any given obstacle function $\psi \in
C(\overline Q)$ such that $\psi \geq c_\psi$ on
$\overline Q$ for a positive constant $c_\psi$
we define a time-dependent constraint 
$K(t)$ by
 $$K(t):=\{z \in V~|~|\nabla z(x)| \le \psi(x,t)
 ~{\rm for ~a.e. ~}x \in \Omega\},~~\forall
 t \in [0,T]. $$
Now, choose a sequence 
$\{\psi_n\}\subset C^1(\overline Q)\}$ so that
 $$ \psi_n \geq c_\psi~{\rm on~} 
 \overline Q,~~\psi_n \to \psi
 ~{\rm in ~}C(\overline Q)~{\rm as~}
 n \to \infty, $$
and consider approximate constraint sets 
 $$K_n(t):=\{z \in V~|~|\nabla z(x)| \le
  \psi_n(x,t)~{\rm for ~a.e. ~}x \in \Omega\},
  ~~\forall t \in [0,T]. $$
Then we have:

(i) $\{K_n(t)\} \in \Psi_S$. The proof is quite
similar to (i) of Example 3.3. We refer to [22;
Lemma 3.1] for the detailed proof.

(ii) $\{K(t)\} \in \Psi_W$. It is enough to check that
$K_n(t)\Longrightarrow K(t)$ on $[0,T]$,
Just as in Example 3.3, this is obtained by using the mapping
 $${\cal F}_\varepsilon z=(1-\varepsilon)z,~~
 \forall z \in V,~~
 \forall~{\rm small ~}\varepsilon >0.$$
 
(iii) Condition (5.3) is satisfied by the strict
positiveness of $\psi$. More precisely, if 
$z \in W_0^{2,q}(\Omega)$, then $|\nabla z| \in 
W_0^{1,q}(\Omega) \subset C(\overline{\Omega})$
and hence there is a positive constant $\kappa$
such that $\kappa |\nabla z| \le c_\psi$ on $Q$
for all $z \in B_W(0)$. Thus (5.3) holds.
\vspace{0.3cm}

\noindent
Therefore, associated with initial datum 
$u_0\in \overline{K(0)}$, the maximal
monotone mapping
$L_{u_0}: D(L_{u_0}) \subset L^2(0,T;V)\to
 L^2(0,T;V^*)$ is well defined on account of
Theorem 3.1.
Next we introduce a semimonotone operator
$A:=A(t,v,u): [0,T]\times H\times V \to V^*$
by:
$$ \langle A(t,v,u), \xi\rangle :=
 \int_\Omega \{a(x,t,v)\nabla u \cdot\nabla \xi
   +b(x,t,v)\xi\}dx $$
  $$\forall v \in H, ~\forall u \in V,~\forall
  \xi \in V,~\forall t \in [0,T], $$
where $a(x,t,v)$ and $b(x,t,v)$ are functions satisfying the
Carath\'eodory condition on $\overline \Omega \times [0,T]\times {\bf R}$
and
 $$ a_*\le a(x,t,v) \le a^*,~-b_* \le b(x,t,v)\le b^*~~{\rm for~a.e.~}(x,t)
  \in \Omega \times (0,T),~{\rm and~all~}v \in {\bf R},$$
for positive constants $a_*,~a^*,~b_*,~b^*$. 

We easily check conditions (a), (b) and (c)
for this operator $A$.
Accordingly we can apply Theorem 5.1
to solve $ f \in L_{u_0}u+ {\cal A}(u,u)$
for given data
$u_0 \in \overline{K(0)}$ 
and $f \in L^2(0,T; V^*)$. This is equivalent
to the weak variational form:
 $$ u \in L^2(0,T;V)\cap C([0,T]; H),~
    |\nabla u|\le \psi~{\rm a.e.~on~}Q,
     ~u(0)=u_0; $$ 
 $$ \int_Q \xi_t(u-\xi)dxdt  
  +\int_Q \{a(x,t,u)\nabla u\cdot 
    \nabla (u-\xi)+b(x,t,u)(u-\xi)\}dxdt $$
  $$ \le \int_0^T \langle f, u-\xi \rangle dt 
  +\frac 12 |u(0)-\xi(0)|^2_{L^2(\Omega)}, $$
  $$ \forall \xi \in L^2(0,T;V)
   \cap W^{1,2}(0,T; H),~|\nabla \eta|
   \le \psi~ {\rm~a.e.~on~}Q.$$
\vspace{0.3cm}

\noindent
{\bf Remark 5.1.} The similar technique in
Application 2 is available for the variational
inequalities arising in models of superconductivity with gradient constraints
or hydrodynamics with velocity
constraints; see [15, 16, 28, 29, 30].
\vspace{1cm}

\begin{center}
{\bf References}
\end{center}

\begin{enumerate}
\item H. W. Alt, An abstract existence theorem for parabolic systems, Comm.
Pure Appl. Anal., 11(2012), 2079-2123.
\item H. Attouch, {\it Variational Convergence for Functions and Operators},
Applicable Math. Ser. Pitman Adv. Pub. Program, Bston-London-Melbourne, 1984.
\item J. P Aubin, Un th\'eor\`eme de 
compacit\'e, C. R. Acad. Sci. Paris,
256 (1963), 5042-5044.
\item V. Barbu, {\it Nonlinear Semi-groupes and
Differential Equations in Banach Spaces},
Noordhoff International Publishing, Leiden, 
1976. 
\item H. Br\'ezis, \'Equations et in\'equations
non lin\'eaires dans les espaces vectoriels en
dualit\'e, Ann. Inst. Fourier, Grenoble
18(1968), 115-175.
\item H. Br\'ezis, Perturbations nonlin\'eaires d'op\'erateurs maximaux monotones, C. R. Acad.,
Sci. Paris, 269(1969), 566-569.
\item H. Br\'ezis, Non linear perturbations of monotone operators, Technical
Report 25, Univ. Kansas, 1972.
\item H. Br\'ezis, Probl\`emes unilat\'eraux, J. Math. pures appl., 51(1972),
1-168.
\item H. Br\'ezis, {\it Op\'eratuers Maximaux Monotones et Semi-groupes de Contractions dans les Espaces de Hilbert}, Math. Studies 5, 
North-Holland, Amsterdam, 1973.
\item F. E. Browder, Nonlinear maximal monotone
operators in Bnach spaces, Math. Ann., 
175(1968), 89-113.
\item F. E. Browder and P. Hess, Nonlinear
mappings of monotone type in Banach spaces, J.
Funct. Anal., 11(1972), 251-294.
\item J. A. Dubinskii, Weak convergence in nonlinear elliptic and parabolic
equations, Amer. Math. Soc. Transl. 2, 67(1968), 226-258.
\item L. C. Evans and R. F. Gariepy, {\it Measure Theory and Fine Properties
of Functions }, CRC Press, Boca Raton-London-New York-Washington, D.C., 1992.
\item A. Friedman, {\it Partial Differential
Equations of Parabolic Type}, Prentice Hall,
1964.
\item T. Fukao and N. Kenmochi, 
Parabolic variational inequalities with weakly time-dependent constraints, 
Adv. Math. Sci. Appl., 23(2013), 365--395. 
\item M. Gokieli, N. Kenmochi and M. 
Niezg\'odka, Variational inequalities of
Navier-Stokes type with time-dependent constraints, J. Math. Anal. Appl. {\bf 449}
(2017), 1229-1247.
\item N. Kenmochi, Nonlinear operators of
 monotone type in reflexive Banach spaces and
nonlinear perturbations, Hiroshima Math. J., 4
(1974), 229-263.
\item N. Kenmochi, R\'esolution de compacit\'e dans les espaces de Banach
d\'ependant du temps, S\'eminaires d'analyse convexe, Montpellier 1979,
Expos\'e 1, 1-26.     .
\item N. Kenmochi, R\'esolution de probl\`emes
variationnels paraboliques non lin\'eaires par
 les m\'ethodes de compacit\'e et monotonie,
Th\`ese de Dcteur de l'Universit\'e, Univ.
Paris VI, 1979.
\item N. Kenmochi, Monotonicity and compactness methods for nonlinear
variational inequalities, pp. 203-298 in {\it Handbook of Differential Equations: Stationary
Partial Differential Equations} {\bf Vol. 4}, Elsevier, Amsterdam, 2007.
\item N. Kenmochi, Solvability of nonlinear evolution equations with time-dependent
constraints and applications, Bull. Fac. Edu., Chiba Univ., 30(1981),1-87.
\item N. Kenmochi, Parabolic quasi-variational
diffusion problems with gradient constraints,
DCDS, Ser. S, {\bf. 6}(2013), 423-438. 
\item N. Kenmochi and M. Niezg\'odka, Weak 
solvability for parabolic variational inclusions
and applications to quasi-variational problems,
Adv. Math. Sci. Appl., 25(2016), 62-97.
\item O. A. Ladyzenskaja, V. A. Solonnikov and
N. N. Ural'ceva, {\it Linear and Quasi-linear
Equations of Parabolic Type}, Transl. 
Mathematical Monographs Vol. 23, Amer. Math. 
Sco., Providence, Rhode Island, 1968.
\item J. Leray and J. L. Lions, Quelques 
r\'esultats de Visik sur les probl\`emes 
elliptiques non lin\'eaires par les m\'ethodes
 de Minty-Browder, Bull. Soc. Math. France
93(1965), 97-107.
\item J. L. Lions, {\it Quelques m\'ethodes de r\'esolution des 
probl\`emes aux 
limites non lin\'eaires}, Dunod Gauthier-Villrs, Paris, 1969.
\item U. Mosco, Convergence of convex sets and of solutions of variational
inequalities, Advances Math. 3(1969), 510-585. 
\item J.-F. Rodrigues and L. Santos, A parabolic
quasi-variational inequality arising in a
superconductivity model, Ann. Scuola Norm. Sup.
Pisa, {\bf 29}(2000), 153-169.
\item L. Santos, A diffusion problem with
gradient constraint and evolutive Dirichlet
condition, Portugal Math., {\bf 48}(1991), 441-
468.
\item L. Santos, Variational problems with
non-constant gradient constraints, Portugal
Math., {\bf 59}(2002), 205-248.
\item J. Simon, Compact sets in the space of $L^p(0,T;B)$, J. Ann.
Mat. pura applic., 146(1986), 65-96.
\end{enumerate}

\end{document}